\documentclass{amsart}
\usepackage{amsmath,amssymb}
\usepackage[margin=1.5cm]{geometry}
\newtheorem{theorem}{Theorem}[section]

\begin{document}
\title[Suzuki polytopes]{On the rank of Suzuki polytopes: an answer to Hubard and Leemans} 

\author[P. Spiga]{Pablo Spiga}
\address{Pablo Spiga,
Dipartimento di Matematica e Applicazioni, University of Milano-Bicocca,\newline
Via Cozzi 55, 20125 Milano, Italy}\email{pablo.spiga@unimib.it}

\begin{abstract}
In this paper we show that the rank of every chiral polytope having a Suzuki group  as automorphism group is $3$. This gives a positive answer to a  conjecture of Isabel Hubard and Dimitri Leemans.
\end{abstract}

\keywords{abstract chiral polytopes, Suzuki simple groups}
\subjclass[2010]{52B11, 20D06}
\maketitle

\section{Introduction}
In~\cite{HL}, Isabel Hubard and Dimitri Leemans have embarked in an extensive analysis on the abstract polytopes admitting an almost simple group with socle a  Suzuki group as automorphism group. Except for an answer to~\cite[Conjecture~1]{HL}, their analysis is very satisfactory and gives a great insight on such abstract polytopes. This his conjecture, Hubard and Leemans ask whether the rank of every chiral polytope having a Suzuki group as automorphism group is $3$. In this paper, we give a positive answer to this conjecture.\footnote{I first heard about this conjecture in the lovely $2014$ summer course in algebraic graph theory by Dimitri Leemans in Rogla, Slovenia. Here I express my gratitude to Dimitri for such a wonderful series of lectures and to the organizers of the Rogla Summer School. } 
\begin{theorem}\label{thrm:main}
Let $\mathrm{Sz}(q)\le G \le \mathrm{Aut}(\mathrm{Sz}(q))$
with $q=2^{2e+1}$, where $e$ is a positive integer and where $\mathrm{Sz}(q)$ is the
 Suzuki group defined over a finite field of cardinality $q$. Then, the rank of any chiral polytope having $G$ as automorphism group is $3$.
\end{theorem}
\section{Notation}

Let $e$ be a positive integer, let  $q:=2^{2e+1}$, let $t:=2^{e}=\sqrt{2q}$, let $\mathbb{F}_q$ be a finite field of cardinality $q$, let $V:=\mathbb{F}_q^4$ be the $4$-dimensional vector space over $\mathbb{F}_4$ consisting of column vectors and let 
$$f:V\times V\to \mathbb{F}_q$$
be the alternating bilinear form whose associated matrix (with respect to the canonical basis of $V$) is 
\[
\iota:=\begin{pmatrix}
0&0&0&1\\
0&0&1&0\\
0&1&0&0\\
1&0&0&0
\end{pmatrix}.
\]

We let $\mathrm{Sp}_4(q)$ be the symplectic group with respect to the bilinear form $f$ and 
 we let $\mathrm{Sz}(q)$ be the Suzuki group, viewed as a subgroup of the symplectic group $\mathrm{Sp}_4(q)$. (Hence, $\mathrm{Sz}(q)$ is the fixed-point subgroup of a graph-field automorphism of $\mathrm{Sp}_4(q)$.)

We now use the characterization of Wilson~\cite{Wilson} of the Suzuki subgroup $\mathrm{Sz}(q)$ in $\mathrm{Sp}_4(q)$. This characterization is  via a suitable commutative product on $V$. Wilson has defined a commutative product
$$\bullet:V\times V\to V$$
by setting
\begin{align*}
e_1\bullet e_1&=0,&e_2\bullet e_1&=e_2,&e_3\bullet e_1&=e_4,&e_4\bullet e_1&=0,\\
e_1\bullet e_2&=e_2,&e_2\bullet e_2&=0,&e_3\bullet e_2&=0,&e_4\bullet e_2&=e_1,\\
e_1\bullet e_3&=e_4,&e_2\bullet e_3&=0,&e_3\bullet e_3&=0,&e_4\bullet e_3&=e_3,\\
e_1\bullet e_4&=0,&e_2\bullet e_4&=e_1,&e_3\bullet e_4&=e_3,&e_4\bullet e_4&=0,\\
\end{align*}
and by extending $\bullet$ to the whole of $V$ using the following formula:
$$\left(\sum_{i=1}^4\lambda_i e_i\right)\bullet\left(\sum_{j=1}^4\mu_je_j\right)=
\sum_{i,j=1}^4\lambda_i^t\mu_j^t (e_i\bullet e_j).$$

Wilson~\cite{Wilson} has shown that $\mathrm{Sz}(q)$ consists exactly of the matrices $g$ in $\mathrm{Sp}_4(q)$ (that is, the matrices  preserving the bilinear form $f$) which satisfy
$$gu\bullet gv=g(u\bullet v), \quad\forall u,v\in V\hbox{ with } f(u,v)=0.$$
(Observe that, there is a typo on~\cite[page~425]{Wilson} and the condition ``$gu\bullet gv=u\bullet v$'' should be replaced by ``$gu\bullet gv=g(u\bullet v)$''. This typo has already be spotted in the mathview by Colva Roney-Dougal.)
We use this characterization of $\mathrm{Sz}(q)$ in the proof of Theorem~\ref{thrm:main}.

\section{Proof of Theorem~$\ref{thrm:main}$}
We argue by contradiction and we suppose that there exists $G$ with $\mathrm{Sz}(q)\le G\le \mathrm{Aut}(\mathrm{Sz}(q))$, for some positive integer $e$, such that $G$ admits a chiral polyhedron of rank different from $3$. In view of~\cite[Theorem~$1$]{HL}, we may suppose that $G=\mathrm{Sz}(q)$ and that $\mathrm{Sz}(q)$ admits a chiral $4$-polytope. In particular, there exist $\sigma_1,\sigma_2,\sigma_3\in \mathrm{Sz}(q)$ with 
$$\mathrm{Sz}(q)=\langle \sigma_1,\sigma_2,\sigma_3\rangle,\quad\hbox{where } \sigma_1\sigma_2\sigma_3,\,\sigma_2\sigma_3,\,\sigma_1\sigma_2\, \hbox{ are involutions of }\mathrm{Sz}(q).$$

Since the involutions in $\mathrm{Sz}(q)$ form a unique conjugacy class~\cite{Suz} and since $\iota\in \mathrm{Sz}(q)$, we may suppose that $$\iota=\sigma_1\sigma_2\sigma_3.$$
Since $\sigma_1\sigma_2=\iota \sigma_3^{-1}$ is an involution, we deduce
$1=(\iota \sigma_3^{-1})^2=\iota\sigma_3^{-1}\iota\sigma_3^{-1}$ and hence $\sigma_3^{-1}\iota \sigma_3^{-1}=\iota$. An entirely similar computation yields
$\sigma_1^{-1}\iota\sigma_1^{-1}=\iota.$
Thus
\begin{equation}\label{eq:2}
\sigma_1^{-1},\sigma_3^{-1}\in \mathcal{X}:=\{x\in \mathrm{Sz}(q)\mid x\iota x=\iota\}.
\end{equation}

In what follows we denote by $x^\perp$ the transpose matrix of $x$. As $\mathrm{Sz}(q)\le \mathrm{Sp}_4(q)$ and as $\iota$ is the matrix associated to the bilinear form $f$, we have
$x^\perp \iota x=\iota$, for every $x\in \mathrm{Sz}(q)$. In particular, for every $x\in\mathcal{X}$, we have $x^\perp \iota x=x\iota x$ and hence
$$x^\perp=x.$$
In particular, $x$ is a symmetric matrix and hence we may write
\begin{align}\label{eq:1}
x&=
\begin{pmatrix}
a_{11}&a_{12}&a_{13}&a_{14}\\
a_{12}&a_{22}&a_{23}&a_{24}\\
a_{13}&a_{23}&a_{33}&a_{34}\\
a_{14}&a_{24}&a_{34}&a_{44}\\
\end{pmatrix}.
\end{align}

We use the characterization of Wilson to determine when a symmetric matrix $x$ as in~\eqref{eq:1} lies in $\mathrm{Sz}(q)$.

From
\begin{align*}
xe_1&=a_{11}e_1+a_{12}e_2+a_{13}e_3+a_{14}e_4,\\
xe_2&=a_{12}e_1+a_{22}e_2+a_{23}e_3+a_{24}e_4,\\
xe_3&=a_{13}e_1+a_{23}e_2+a_{33}e_3+a_{34}e_4,\\
xe_4&=a_{14}e_1+a_{24}e_2+a_{34}e_3+a_{44}e_4,
\end{align*}
we deduce
\begin{align*}
xe_2=&x(e_1\bullet e_2)=xe_1\bullet xe_2
=(a_{11}e_1+a_{12}e_2+a_{13}e_3+a_{14}e_4)\bullet (a_{12}e_1+a_{22}e_2+a_{23}e_3+a_{24}e_4),\\
=&
a_{11}^ta_{12}^t(e_1\bullet e_1)+a_{11}^ta_{22}^t(e_1\bullet e_2)+a_{11}^ta_{23}^t(e_1\bullet e_3)+a_{11}^ta_{24}^t(e_1\bullet e_4)\\
&+a_{12}^ta_{12}^t(e_2\bullet e_1)+a_{12}^ta_{22}^t(e_2\bullet e_2)+a_{12}^ta_{23}^t(e_2\bullet e_3)+a_{12}^ta_{24}^t(e_2\bullet e_4)\\
&+a_{13}^ta_{12}^t(e_3\bullet e_1)+a_{13}^ta_{22}^t(e_3\bullet e_2)+a_{13}^ta_{23}^t(e_3\bullet e_3)+a_{13}^ta_{24}^t(e_3\bullet e_4)\\
&+a_{14}^ta_{12}^t(e_4\bullet e_1)+a_{14}^ta_{22}^t(e_4\bullet e_2)+a_{14}^ta_{23}^t(e_4\bullet e_3)+a_{14}^ta_{24}^t(e_4\bullet e_4)\\
=&
a_{11}^ta_{22}^t(e_1\bullet e_2)+a_{11}^ta_{23}^t(e_1\bullet e_3)
+a_{12}^ta_{12}^t(e_2\bullet e_1)+a_{12}^ta_{24}^t(e_2\bullet e_4)\\
&+a_{13}^ta_{12}^t(e_3\bullet e_1)+a_{13}^ta_{24}^t(e_3\bullet e_4)+
a_{14}^ta_{22}^t(e_4\bullet e_2)+a_{14}^ta_{23}^t(e_4\bullet e_3)\\
=&a_{11}^ta_{22}^t e_2+a_{11}^ta_{23}^t e_4
+a_{12}^ta_{12}^t e_2+a_{12}^ta_{24}^t e_1
+a_{13}^ta_{12}^t e_4+a_{13}^ta_{24}^t e_3
+a_{14}^ta_{22}^t e_1+a_{14}^ta_{23}^t e_3\\
=&
(a_{12}^ta_{24}^t+a_{14}^ta_{22}^t) e_1+
(a_{11}^ta_{22}^t +a_{12}^ta_{12}^t)e_2+(a_{13}^ta_{24}^t
+a_{14}^ta_{23}^t )e_3+(a_{11}^ta_{23}^t+a_{13}^ta_{12}^t )e_4\\
\end{align*}
Therefore we obtain the system of equations:
\begin{align}\label{eq:31}
a_{12}^ta_{24}^t+a_{14}^ta_{22}^t&=a_{12},\\\label{eq:32}
a_{11}^ta_{22}^t +a_{12}^{2t}&=a_{22},\\\label{eq:33}
a_{13}^ta_{24}^t+a_{14}^ta_{23}^t &=a_{23},\\\label{eq:34}
a_{11}^ta_{23}^t+a_{13}^ta_{12}^t &=a_{24}
\end{align}
An entirely analogous computation can be done imposing 
\begin{align*}
&xe_1\bullet xe_3=x(e_1\bullet e_3)=xe_4,\,\,\,\, xe_4\bullet xe_3=x(e_4\bullet e_3)=xe_3,\,\,\,\, xe_1\bullet xe_3=x(e_2\bullet e_4)=xe_1,\\
&xe_1\bullet xe_4=x(e_1\bullet e_4)=x\cdot 0=0,\,\,\,\, xe_2\bullet xe_3=x(e_2\bullet e_3)=x\cdot 0=0.
\end{align*}
We skip the computations above and we report the resulting equations:
\begin{align}\label{eq:35}
a_{12}^ta_{34}^t+a_{14}^ta_{23}^t&=a_{14},\\\label{eq:36}
a_{13}^ta_{34}^t+a_{14}^ta_{33}^t &=a_{34},\\\label{eq:38}
a_{11}^ta_{33}^t+a_{13}^{2t} &=a_{44}\\\label{eq:39}
a_{22}^ta_{44}^t+a_{24}^{2t}&=a_{11},\\\label{eq:310}
a_{23}^ta_{44}^t+a_{24}^ta_{34}^t &=a_{13},\\\label{eq:312}
a_{33}^ta_{44}^t+a_{34}^{2t} &=a_{33},\\\label{eq:313}
a_{14}a_{24}&=a_{12}a_{44},\\\label{eq:314}
a_{12}a_{14}&=a_{11}a_{24},\\\label{eq:315}
a_{13}a_{44}&=a_{14}a_{34},\\\label{eq:316}
a_{11}a_{34}&=a_{13}a_{14},\\\label{eq:317}
a_{22}a_{34}&=a_{24}a_{23},\\\label{eq:318}
a_{22}a_{13}&=a_{12}a_{23},\\\label{eq:319}
a_{23}a_{34}&=a_{24}a_{33},\\\label{eq:320}
a_{12}a_{33}&=a_{13}a_{23}.
\end{align}

Writing the equality $x\iota x=\iota$ explicitly, we obtain
\begin{align*}
\begin{pmatrix}
0&0&0&1\\
0&0&1&0\\
0&1&0&0\\
1&0&0&0
\end{pmatrix}&=\begin{pmatrix}
a_{11}&a_{12}&a_{13}&a_{14}\\
a_{12}&a_{22}&a_{23}&a_{24}\\
a_{13}&a_{23}&a_{33}&a_{34}\\
a_{14}&a_{24}&a_{34}&a_{44}
\end{pmatrix}
\begin{pmatrix}
0&0&0&1\\
0&0&1&0\\
0&1&0&0\\
1&0&0&0
\end{pmatrix}
\begin{pmatrix}
a_{11}&a_{12}&a_{13}&a_{14}\\
a_{12}&a_{22}&a_{23}&a_{24}\\
a_{13}&a_{23}&a_{33}&a_{34}\\
a_{14}&a_{24}&a_{34}&a_{44}
\end{pmatrix}\\
&=
\begin{pmatrix}
a_{14}&a_{13}&a_{12}&a_{11}\\
a_{24}&a_{23}&a_{22}&a_{12}\\
a_{34}&a_{33}&a_{23}&a_{13}\\
a_{44}&a_{34}&a_{24}&a_{14}
\end{pmatrix}
\begin{pmatrix}
a_{11}&a_{12}&a_{13}&a_{14}\\
a_{12}&a_{22}&a_{23}&a_{24}\\
a_{13}&a_{23}&a_{33}&a_{34}\\
a_{14}&a_{24}&a_{34}&a_{44}
\end{pmatrix}.
\end{align*}
By computing explicitly the matrix product on the right hand side and by comparing the coefficients in position $(1,4)$ and $(2,3)$ with the matrix on the left hand side, we obtain two further equations
\begin{align}\label{eq:321}
a_{14}^2+a_{13}a_{24}+a_{12}a_{34}+a_{11}a_{44}&=1,\\\label{eq:322}
a_{24}a_{13}+a_{23}^2+a_{22}a_{33}+a_{12}a_{34}&=1.
\end{align}

We now start solving this system of equations by dividing our task in various cases.

\subsection{Case $a_{24}\ne 0$} From~\eqref{eq:313},~\eqref{eq:314},~\eqref{eq:317} and~\eqref{eq:319}, we deduce
\begin{equation}\label{new}
a_{14}=\frac{a_{12}a_{44}}{a_{24}},\,\,a_{11}=\frac{a_{12}^2a_{44}}{a_{24}^2},\,\,a_{23}=\frac{a_{22}a_{34}}{a_{24}},\,\,a_{33}=\frac{a_{22}a_{34}^2}{a_{24}^2}.
\end{equation}
Using this value of $a_{14}$ and $a_{11}$ in~\eqref{eq:321}, we obtain
\begin{equation}\label{eq:323}a_{13}a_{24}+a_{12}a_{34}=1.
\end{equation}
From~\eqref{eq:315}, we get
$$a_{13}a_{44}=\frac{a_{12}a_{44}a_{34}}{a_{24}}.$$
Using this equality and multiplying~\eqref{eq:323} by $a_{44}/a_{24}$, we deduce
$$\frac{a_{44}}{a_{24}}=\frac{a_{44}}{a_{24}}\left(a_{13}a_{24}+a_{12}a_{34}\right)=a_{13}a_{44}+\frac{a_{12}a_{44}a_{34}}{a_{24}}=0.$$
Therefore $a_{44}=0$. Using this value of $a_{44}$ in~\eqref{new}, we deduce $a_{14}=a_{11}=0$.

Substituting the value of $a_{23}$ from~\eqref{new} in~\eqref{eq:318}, we get $a_{22}a_{13}=\frac{a_{12}a_{22}a_{34}}{a_{24}}$. Therefore, multiplying~\eqref{eq:323} by $a_{22}/a_{24}$, we obtain 
$$\frac{a_{22}}{a_{24}}=\frac{a_{22}}{a_{24}}\left(a_{13}a_{24}+a_{12}a_{34}\right)=a_{22}a_{13}+\frac{a_{12}a_{22}a_{34}}{a_{24}}=0.$$
Therefore $a_{22}=0$. Using this value of $a_{22}$ in~\eqref{eq:32}, we deduce $a_{12}=0$. Now,~\eqref{new} gives $a_{33}=0$. From~\eqref{eq:38}, we get $a_{13}=0$. Summing up,
$$a_{11}=a_{12}=a_{13}=a_{14}=0,$$
contradicting the fact that $x$ is a non-singular matrix.

\subsection{Case $a_{24}= 0$ and $a_{12}\ne 0$} From~\eqref{eq:313},~\eqref{eq:314}, we obtain
$$a_{44}=0,\,a_{14}=0.$$
Now,~\eqref{eq:321} gives $a_{12}a_{34}=1$. However,~\eqref{eq:31} with $a_{24}=a_{14}=0$ gives $a_{12}=0$, which is a contradiction.

\subsection{Case $a_{24}=a_{12}= 0$ and $a_{34}\ne 0$} From~\eqref{eq:317}, we have
$a_{22}=0$. Using this, from~\eqref{eq:322}, we get $a_{23}=1$. However,~\eqref{eq:319} yields $a_{23}a_{34}=0$, which is a contradiction.

\subsection{Case $a_{24}=a_{12}=a_{34}= 0$ and $a_{13}\ne 0$} From~\eqref{eq:315}, we have $a_{44}=0$. Now,~\eqref{eq:321} yields $a_{14}=1$. However,~\eqref{eq:316} gives $a_{13}a_{14}=0$, which is a contradiction.

\subsection{Case $a_{24}=a_{12}=a_{34}=a_{13}= 0$ and $a_{14}\ne 0$}From~\eqref{eq:31}, we have $a_{22}=0$. Now,~\eqref{eq:322} gives $a_{23}=1$. In turn, from~\eqref{eq:34}, we obtain $a_{11}=0$. Using this, from~\eqref{eq:321}, we get $a_{14}=1$. From~\eqref{eq:36} and~\eqref{eq:38}, we get $a_{33}=a_{44}=0$. Therefore
$$x=\iota$$
and we have our first element in $\mathcal{X}$.

\subsection{Case $a_{24}=a_{12}=a_{34}=a_{13}= a_{14}= 0$}From~\eqref{eq:321}, we have $a_{44}=a_{11}^{-1}$. Now,~\eqref{eq:33} yields $a_{23}=0$. In turn, from~\eqref{eq:322}, we have $a_{33}=a_{22}^{-1}$. Finally,~\eqref{eq:39} yields $a_{22}^ta_{11}^{-t}=a_11$, that is, $a_{22}^t=a_{11}^{t+1}$. Raising both sides of this equality by $2t$ and recallying that $2t^2=q$, we obtain $a_{22}=a_{11}^{2t+1}$. Therefore
\[
x=\begin{pmatrix}
a_{11}&0&0&0\\
0&a_{11}^{2t+1}&0&0\\
0&0&a_{11}^{-2t-1}&0\\
0&0&0&a_{11}^{-1}
\end{pmatrix}.
\]

Summing up, we have shown that
$$\mathcal{X}=\left\{
\begin{pmatrix}
0&0&0&1\\
0&0&1&0\\
0&1&0&0\\
1&0&0&0
\end{pmatrix}
\right\}\cup\left\{
\begin{pmatrix}
a&0&0&0\\
0&a^{2t+1}&0&0\\
0&0&a^{-2t-1}&0\\
0&0&0&a^{-1}
\end{pmatrix}\mid a\in \mathbb{F}_q\setminus\{0\}
\right\}.$$
\subsection{Pulling the threads of the argument}We now use~\eqref{eq:2}. If $\sigma_1^{-1}$ or $\sigma_3^{-1}$ equals $\iota$, then we obtain that either $\sigma_1\sigma_2=1$ or $\sigma_2\sigma_3=1$, contradicting the fact that $\sigma_1\sigma_2$ and $\sigma_2\sigma_3$ are both involutions. Therefore,
\begin{align*}
\sigma_1&=\begin{pmatrix}
a^{-1}&0&0&0\\
0&a^{-2t-1}&0&0\\
0&0&a^{2t+1}&0\\
0&0&0&a
\end{pmatrix}\\
\sigma_3&=\begin{pmatrix}
b^{-1}&0&0&0\\
0&b^{-2t-1}&0&0\\
0&0&b^{2t+1}&0\\
0&0&0&b
\end{pmatrix},
\end{align*}
for some $a,b\in\mathbb{F}_q\setminus\{0\}$. Observe that $\sigma_1$ and $\sigma_3$ centralize each other. In particular, $\langle \sigma_1,\sigma_3\rangle$ is abelian. Now, $$\mathrm{Sz}(q)=\langle\sigma_1,\sigma_2,\sigma_3\rangle=\langle \sigma_1,\sigma_3,\sigma_1\sigma_2\sigma_3\rangle=\langle\sigma_1,\sigma_3,\iota\rangle.$$  
An elementary computation gives $\sigma_1^\iota=\sigma_1^{-1}$ and $\sigma_3^{\iota}=\sigma_3^{-1}$. Therefore, $\langle \sigma_1,\iota,\sigma_3\rangle$ is soluble. As $\mathrm{Sz}(q)$ is a non-abelian simple group, we obtain our final contradiction.


\thebibliography{10}
\bibitem{HL}I.~Hubard, D.~Leemans, Chiral polytopes and Suzuki simple groups. Rigidity and symmetry, 155--175, \textit{Fields Inst. Commun.} \textbf{70}, Springer, New York, 2014. 
\bibitem{Suz}M.~Suzuki, On a class of doubly transitive groups, \textit{Ann. Math.}  \textbf{75}, 105--145, 1962.
\bibitem{Wilson}R.~A.~Wilson, A new approach to the Suzuki groups, \textit{Math. Proc. Cambridge Philos. Soc.} \textbf{148} (2010), 425--428.
\end{document}